\documentclass[12pt]{article}
\usepackage[dvipsnames]{xcolor}
\usepackage{amsmath}
\usepackage{setspace}
\usepackage{bm}
\usepackage{bbm}
\usepackage{amssymb}
\usepackage{indentfirst}
\usepackage[superscript]{cite}
\usepackage{geometry}
\usepackage{mathtools}
\usepackage{hyperref}

\newcommand*{\citena}[1]{%
\begingroup
[\color{Green}
\romannumeral-`\x 
\setcitestyle{numbers}%
\cite{#1}%
\endgroup
]\ignorespacesafterend
}

\newcommand*{\citesup}[1]{%
\begingroup
\color{Green}
\cite{#1}%
\endgroup
\ignorespacesafterend
}

\newcommand*{\eqrefe}[1]{%
\begingroup
(\color{BrickRed}
\romannumeral-`\x 
\setcitestyle{numbers}%
\ref{eq:#1}%
\endgroup
)\ignorespacesafterend
}

\newcommand*{\secrefe}[1]{%
\begingroup
(\color{Aquamarine}
\romannumeral-`\x 
\setcitestyle{numbers}%
\ref{#1}%
\endgroup
)\ignorespacesafterend
}

\newtagform{Tags}[\textcolor{BrickRed}]{\color{Black}(}{)}

\usepackage[super]{natbib}
\setcitestyle{super}

\newcommand{\ii}{\bm{i}}

\DeclareMathOperator{\csch}{csch}

\begin{document}
\title{Integrals for the Zeta Function and its Generalizations}
\date{June 22, 2022}
\author{
José Risomar Sousa\\
\href{https://orcid.org/0009-0002-7779-0808}
     {ORCID: 0009-0002-7779-0808}
}
\maketitle
\usetagform{Tags}

\begin{abstract}
A formula for the Riemann zeta function is obtained as an extension of the Faulhaber formula and, from there, formulae for its generalizations (the Hurwitz zeta, $\zeta(-k,b)$, the polylogarithm, $\mathrm{Li}_{-k}(e^z)$, and the Lerch transcendent, $\Phi(e^z,-k,b)$), are derived that coincide with their Abel-Plana expressions. It allows one to find, for example, the Taylor series expansion of $H_{-k}(n)$ about $n=0$ (when $k$ is a positive integer, a finite Taylor series is obtained, which is nothing but the Faulhaber formula). The used method requires evaluating the limit of $\Phi\left(e^{2\pi\ii\,x},-2k+1,n+1\right)+\pi\ii\,x\,\Phi\left(e^{2\pi\ii\,x},-2k,n+1\right)/k$ when $x$ goes to zero, which is an interesting problem in itself.
\end{abstract}

\tableofcontents

\section{Introduction}
The formula derivations that follow next are based on two main ideas that were introduced in two previous papers, \citena{Faulhaber} and \citena{Lerch_neg}. In this paper, those ideas are explored further to see what new results can be obtained. The reference \citena{Butzer} by Butzer et al. (2011) is highly relevant, since it discusses summation formulae (Euler–Maclaurin, Abel–Plana, Poisson) and their connections with integral transforms.\\

The Faulhaber formula is a closed formula for the harmonic numbers, $H_{-k}(n)$, when $k$ is a positive integer. It is the summation of positive integer powers of consecutive integers starting at one, and can be expressed as a function of Bernoulli numbers --- the numbers that appear in the Taylor series expansion of $x/(e^x-1)$. For positive integer odd powers, this expression is:
\begin{equation} \nonumber
H_{-2k+1}(n)=\sum_{j=1}^{n}j^{2k-1}=\frac{n^{2k-1}}{2}+(2k-1)!\sum_{j=0}^{k-1}\frac{B_{2j}n^{2k-2j}}{(2j)!(2k-2j)!}
\end{equation}
\indent As the reader may know, when $\Re(k)<1$, the limit of $H_{-k}(n)$ as $n$ approaches infinity is the Riemann zeta function, so these two functions are closely interconnected.\\

The first of the aforementioned ideas is to use the analytic continuation of the Bernoulli numbers, achievable through the zeta function, to extend the Faulhaber formula. Since,
\begin{equation} \nonumber
\frac{B_{2j}}{(2j)!}=-2(-1)^j(2\pi)^{-2j}\zeta(2j) \text{,}
\end{equation}
\noindent the modified form below can be obtained,
\begin{equation} \nonumber
H_{-2k+1}(n)=\frac{n^{2k}}{2k}+\frac{n^{2k-1}}{2}+2(2k-1)!(2\pi\ii)^{-2k}\zeta(2k)-2(2k-1)!\,n^{2k}\sum_{j=1}^{k}\frac{(2\pi\ii\,n)^{-2j}\zeta(2j)}{(2k-2j)!}
\end{equation}
\indent The whole focus is then on how one can find ways to obtain a closed form for the key summation below:
\begin{equation} \label{eq:key_sum} 
\sum_{j=1}^{k}\frac{(2\pi\ii\,n)^{-2j}\zeta(2j)}{(2k-2j)!}
\end{equation}
\indent The second idea, first introduced in \citena{Lerch_neg}, is the following straightforward identity, which only holds for the analytic continuation of the Lerch $\Phi$ function at the negative integers $-k$:
\begin{equation} \label{eq:Soma_Lerch}
\Phi(e^{z},-k,u+v)=k!\,\sum_{j=0}^{k}\frac{\Phi\left(e^{z},-j,v\right)u^{k-j}}{j!(k-j)!}
\end{equation}
\indent A special case of this identity is the summation of polylogs, when $v=1$:
\begin{equation} \label{eq:Soma_Polylog}
e^{z}\,\Phi(e^{z},-k,b+1)=k!\,\sum_{j=0}^{k}\frac{\mathrm{Li}_{-j}\left(e^{z}\right)b^{k-j}}{j!(k-j)!}
\end{equation}
\indent Since $e^z=1$ leads to singularities, the analogous expression for  the Hurwitz zeta function at the negative integers was made possible through some relations available in the literature, as explained in \citena{Lerch_neg}, and is slightly different:
\begin{equation} \label{eq:Soma_Hurwitz} 
\zeta(-k,u+v)=-\frac{u^{k+1}}{k+1}+k!\,\sum_{j=0}^{k}\frac{\zeta(-j,v)\,u^{k-j}}{j!(k-j)!} \text{,} 
\end{equation}

\section{Riemann zeta function} \label{preface}
Leveraging the two aforementioned ideas, one starts with an optimal expression for the Riemann zeta function and then from it the Hurwitz zeta, $\zeta(-k,b)$, and the polylogarithm, $\mathrm{Li}_{-k}(e^z)$, are obtained. Lastly, from the polylogarithm, the most convoluted one, the Lerch $\Phi$ function\citesup{Lerch1887}, $\Phi(e^z,-k,b)$, is obtained.

\subsection{Integral from the literature} \label{Faulhaber2}
Let us see how a different expression for $H_{-k}(n)$ can be derived than the one that was created previously in \citena{Faulhaber}, by using a different integral for the zeta function.\\

First, a closed form for \eqrefe{key_sum} is needed. If $\Re{(k)}>1$, one integral representation of the zeta function available in the literature\citesup{Abramowitz} is:
\begin{equation} \label{eq:zeta(k)pos}
\zeta(k)=\frac{1}{(k-1)!}\int_{0}^{\infty}\frac{x^{k-1}}{e^x-1}\,dx=\frac{1}{(k-1)!}\int_{0}^{1}\frac{\left(-\log{u}\right)^{k-1}}{1-u}\,du \text{}
\end{equation}
\indent Hence, using equation \eqrefe{zeta(k)pos}, one has:
\begin{equation} \nonumber
\sum_{j=1}^{k}\frac{(2\pi\ii\,n)^{-2j}\zeta(2j)}{(2k-2j)!}=\int_{0}^{1}\sum_{j=1}^{k}\frac{(2\pi\ii\,n)^{-2j}}{(2k-2j)!}\frac{\left(-\log{u}\right)^{2j-1}}{(2j-1)!(1-u)}\,du
\end{equation}
\indent By replacing the summation within the above integral with a closed form its extension is obtained, which this time is slightly simpler than before:
\begin{multline} \nonumber
\sum_{j=1}^{k}\frac{(2\pi\ii\,n)^{-2j}\zeta(2j)}{(2k-2j)!}=-\frac{(2\pi\ii\,n)^{-2k}}{2(2k-1)!}\int_{0}^{1}\frac{\left(2\pi\ii\,n+\log{u}\right)^{2k-1}+\left(-2\pi\ii\,n+\log{u}\right)^{2k-1}}{1-u}\,du\\=\frac{(2\pi\ii\,n)^{-2k}}{2(2k-1)!}\int_{0}^{\infty}\frac{\left(2\pi\ii\,n+x\right)^{2k-1}+\left(-2\pi\ii\,n+x\right)^{2k-1}}{e^x-1}\,dx
\end{multline}
\indent Now, repeating the same process outlined in reference \citena{Faulhaber}, and making a change of variables, one obtains,
\begin{equation} \label{eq:HN_final_v2}
\sum_{j=1}^{n}j^{k}=\frac{n^{k+1}}{k+1}+\frac{n^{k}}{2}+\zeta(-k)-n^{k+1}\int_{0}^{\pi/2}\left(1-\coth{\left(\pi\,n\tan{v}\right)}\right)\frac{\sin{k\,v}}{(\cos{v})^{k+2}}\,dv \text{,}
\end{equation}
\noindent which is a bit simpler than the previous formula:
\begin{equation} \label{eq:HN_final_v1}
\sum_{j=1}^{n}j^{k}=\frac{n^{k+1}}{k+1}+\frac{n^{k}}{2}+\zeta(-k)+\frac{\pi\,n^{k+2}}{k+1}\int_{0}^{\pi/2}\left(\sec{v}\csch{\left(\pi\,n\tan{v}\right)}\right)^2\left(1-\frac{\cos{(k+1)\,v}}{(\cos{v})^{k+1}}\right)\,dv \text{}
\end{equation}
\indent As is plain to see, these formulae are by no means unique. None of these two forms allows the Taylor series expansion of $H_{-k}(n)$ about $n=0$ to be calculated, if $\Re(k)<0$. They do not produce the right analytic continuation of $H_{-k}(n)$ over $n$ either, if $\Re{(n)} \le 0$ (but over $k$ it seems that they always do).

\subsection{A simpler expression for $\zeta(k)$} \label{Simpler_zeta}
The integral representations that were used for the zeta function in \secrefe{Faulhaber2} tend to make the final formulae a bit complicated, so the objective is to try and see if a simpler expression for \eqrefe{key_sum} can be found. But in order to do that, first it is necessary to find a simpler expression for $\zeta(k)$.\\

From the literature, a Taylor series expansion for $x\cot{x}$\citesup{Abramowitz} is known, whose $k$-th derivatives give the Bernoulli numbers, and hence also the zeta function at the even integers:
\begin{equation} \nonumber
x\cot{x}=\sum_{k=0}^{\infty}\frac{(-1)^k B_{2k}(2x)^{2k}}{(2k)!} \text{, if } |x|<\pi
\end{equation}
\indent But one does not need to use this function, a better option is the generating function of the zeta function at the even integers:
\begin{equation} \label{eq:Even_Zeta_GF}
\sum_{k=0}^{\infty}\zeta(2k)x^{2k}=-\frac{\pi x\cot{\pi x}}{2}
\end{equation}

\subsection{Derivatives of trigonometric functions} \label{Trigon}
\indent In \citena{Lerch_neg}, an expression for the $k$-th derivatives of the cotangent was created, which proves to be very useful now. The formula below makes use of the Kronecker delta function:
\begin{equation} \label{eq:Cot_final}
\frac{\text{d}^k(\cot{a\,x})}{\text{d}\,x^k}=-\ii\,\delta_{0\,k}-2\,\ii(2\,\ii\,a)^k\,\mathrm{Li}_{-k}\left(e^{2\,\ii\,a\,x}\right) \text{, where } \delta_{0\,k}=1 \text{ iff } k=0 
\end{equation}
\indent For the record, the expressions for the derivatives of the other trigonometric functions are (the translated arc formulae allow one type to be converted into another):
\begin{equation} \label{eq:Cosec_final} \nonumber
\frac{\text{d}^k}{\text{d}\,x^k}\left(\frac{1}{\sin{a\,x}}\right)=-2\,\ii(2\,\ii\,a)^k\,e^{\ii\,a\,x}\,\Phi\left(e^{2\,\ii\,a\,x},-k,\frac{1}{2}\right) \text{}
\end{equation}
\begin{equation} \label{Tan_final} \nonumber
\frac{\text{d}^k\left(\tan{(a\,x+b)}\right)}{\text{d}\,x^k}=\ii\,\delta_{0\,k}+2\,\ii(2\,\ii\,a)^k\,\mathrm{Li}_{-k}\left(-e^{2\,\ii(a\,x+b)}\right) 
\end{equation}
\begin{equation} \label{eq:Sec_final} \nonumber
\frac{\text{d}^k}{\text{d}\,x^k}\left(\frac{1}{\cos{(a\,x+b)}}\right)=2\,(2\,\ii\,a)^k\,e^{\ii(a\,x+b)}\,\Phi\left(-e^{2\,\ii(a\,x+b)},-k,\frac{1}{2}\right) \text{}
\end{equation}\\
\indent The problem with equation \eqrefe{Cot_final} is the fact that it is improper at $x=0$, but from \eqrefe{Even_Zeta_GF} we know the limits exist. So one needs to differentiate $x\cot{\pi x}$ with respect to $x$ and take the limit as $x$ goes to zero.\\

The Leibniz rule states that the $j$-th derivative of $x f(x)$ is $j f^{(j-1)}(x)+x f^{(j)}(x)$. Therefore, for $j$ a non-negative integer, the zeta function at the even integers can be written as:
\begin{multline} \label{eq:zeta_at_even_int}
\lim_{x\to 0}-\frac{\pi}{2\,j!}\left(-j\left(\ii\,\delta_{0\,j-1}+2\,\ii(2\,\pi\ii)^{j-1}\,\mathrm{Li}_{-j+1}\left(e^{2\,\pi\ii\,x}\right)\right)-x\left(\ii\,\delta_{0\,j}+2\,\ii(2\,\pi\ii)^j\,\mathrm{Li}_{-j}\left(e^{2\,\pi\ii\,x}\right)\right)\right)\\=\begin{cases} \nonumber 
\zeta(j), & \text{if $j$ is even}\\
0, & \text{if $j$ is odd} 
\end{cases}
\end{multline}
\indent Therefore, the key summation from equation \eqrefe{key_sum} that must be evaluated is:\footnote{Note the indicator function is being used: $\mathbbm{1}_{2|j}=1$ if $j$ is even, 0 otherwise.}
\begin{multline} \nonumber
\sum_{j=0}^{k}\frac{(2\pi\ii\,n)^{-2j}\zeta(2j)}{(2k-2j)!}=\sum_{j=0}^{2k}\frac{\mathbbm{1}_{2|j}\cdot(2\pi\ii\,n)^{-j}\zeta(j)}{(2k-j)!}\\=\frac{\pi}{2}\lim_{x\to 0}\sum_{j=0}^{2k}\frac{(2\pi\ii\,n)^{-j}}{(2k-j)!j!}\left(j\left(\ii\,\delta_{0\,j-1}+2\ii(2\pi\ii)^{j-1}\,\mathrm{Li}_{-j+1}\left(e^{2\pi\ii\,x}\right)\right)+x\left(\ii\,\delta_{0\,j}+2\ii(2\pi\ii)^j\,\mathrm{Li}_{-j}\left(e^{2\pi\ii\,x}\right)\right)\right)
\end{multline}
\indent Hence, the limit of the expression below must be evaluated:
\begin{equation} \nonumber
\lim_{x\to 0}\frac{1}{4n(2k-1)!}+\frac{\pi\ii\,x}{2(2k)!}+\pi\ii\sum_{j=0}^{2k}\frac{(2\pi\ii\,n)^{-j}}{(2k-j)!j!}\left(j\,(2\pi\ii)^{j-1}\,\mathrm{Li}_{-j+1}\left(e^{2\pi\ii\,x}\right)+x\,(2\pi\ii)^j\,\mathrm{Li}_{-j}\left(e^{2\pi\ii\,x}\right)\right)
\end{equation}
\indent That expression can be simplified as below:
\begin{equation} \label{eq:polylog_sum}
\lim_{x\to 0}\frac{1}{4n(2k-1)!}+\frac{1}{2}\sum_{j=1}^{2k}\frac{n^{-j}\,\mathrm{Li}_{-j+1}\left(e^{2\pi\ii\,x}\right)}{(2k-j)!(j-1)!}+\pi\ii\,x\sum_{j=0}^{2k}\frac{n^{-j}\,\mathrm{Li}_{-j}\left(e^{2\pi\ii\,x}\right)}{(2k-j)!j!}
\end{equation}

\subsection{A limit involving the Lerch $\Phi$} \label{Lerch_diff_lim_sec}
The evaluation of this limit is related to the analytic continuation of summation formulae, as discussed in the context of Poisson and Abel–Plana summation in \citena{Butzer}.\\

By applying the identity \eqrefe{Soma_Polylog} to \eqrefe{polylog_sum} (disregarding the trivial term), the limit reduces to the following expression,
\begin{equation} \nonumber
\lim_{x\to 0}\frac{n^{-2k}e^{2\pi\ii\,x}}{2(2k-1)!}\,\Phi\left(e^{2\pi\ii\,x},-2k+1,n+1\right)+\pi\ii\,x\frac{n^{-2k}e^{2\pi\ii\,x}}{(2k)!}\,\Phi\left(e^{2\pi\ii\,x},-2k,n+1\right) \text{,}
\end{equation}
\noindent which can be rewritten as:
\begin{equation} \label{eq:Lerch_diff_lim}
\lim_{x\to 0}\frac{n^{-2k}}{2(2k-1)!}\left(\Phi\left(e^{2\pi\ii\,x},-2k+1,n+1\right)+\frac{\pi\ii\,x}{k}\Phi\left(e^{2\pi\ii\,x},-2k,n+1\right)\right)
\end{equation}\\
\indent At this point, one might think it could help to use the closed form of the Lerch $\Phi$ at the negative integers from reference \citena{Lerch_neg},
\begin{equation} \label{eq:Lerch__closed_form}
\Phi(e^{2\pi\ii\,x},-k,\,n+1)=-\frac{1}{e^{2\pi\ii\,x}-1}\sum_{q=0}^{k}\left(\frac{e^{2\pi\ii\,x}}{e^{2\pi\ii\,x}-1}\right)^q\sum_{j=0}^{q}\binom{q}{j}(-1)^j (j+n+1)^{k} \text{,}
\end{equation}
\noindent but it is in fact very cumbersome in this case. It is not hard to approximate $1/(e^{2\pi\ii\,x}-1)$ or $1/(e^{-2\pi\ii\,x}-1)$ when $x$ is small, but their powers have patterns that are hard to find (although not impossible).\\

The right way to solve this problem is to use an alternative integral for the Lerch $\Phi$ that holds at the negative integers, whose derivation is explained in section \secrefe{Lerch_used}. The integral below holds for every integer $k$:
\begin{multline} \label{eq:Lerch_final2} \nonumber
\Phi(e^{2\pi\ii\,x},-k,\,n+1)=\frac{(n+1)^k}{2}+(-2\pi\ii\,x)^{-k-1}\,e^{-2\pi\ii\,x(n+1)}\,\Gamma\left(k+1,-2\pi\ii\,x(n+1)\right)\\
+\int_{0}^{\pi/2}\frac{1-\coth{\left(\pi\,n\tan{v}\right)}}{(\cos{v})^{2}}\left((n+1)^2+(\tan{v})^2\right)^{k/2}\sin{\left(k\arctan{\frac{\tan{v}}{n+1}}+2\pi\ii\,x\tan{v}\right)}\,dv
\end{multline}\\
\indent After some analysis, one finds that the non-trivial term of the limit in \eqrefe{Lerch_diff_lim} can be attributed to the upper incomplete gamma function in this formula. That is, for the other terms, one can simply evaluate the expression at $x=0$ (it is not improper anymore).\\

Evaluated at $x=0$, the integral that does not vanish becomes:
\begin{multline} \nonumber
\int_{0}^{\pi/2}\frac{1-\coth{\left(\pi\tan{v}\right)}}{(\cos{v})^{2}}\left((n+1)^2+(\tan{v})^2\right)^{(2k-1)/2}\sin{\left((2k-1)\arctan{\frac{\tan{v}}{n+1}}\right)}\,dv\\
=-\frac{\ii}{2}\int_{0}^{\pi/2}\frac{1-\coth{\left(\pi\tan{v}\right)}}{(\cos{v})^{2}}\left((n+1+\ii\tan{v})^{2k-1}-(n+1-\ii\tan{v})^{2k-1}\right)\,dv\\
=-\frac{\ii}{2}\int_{0}^{\infty}(1-\coth{\pi x})\left((n+1+\ii\,x)^{2k-1}-(n+1-\ii\,x)^{2k-1}\right)\,dx
\end{multline}
\indent For the limit of the difference between the upper incomplete gamma functions one has:
\begin{equation} \nonumber
\lim_{x\to 0}\frac{e^{-2\pi\ii\,x(n+1)}}{(2\pi\ii\,x)^{2k}}\left(\Gamma\left(2k,-2\pi\ii\,x(n+1)\right)-\frac{1}{2\,k}\Gamma\left(2k+1,-2\pi\ii\,x(n+1)\right)\right)=-\frac{(n+1)^{2k}}{2k}
\end{equation}\\
\indent Therefore, the limit that one is looking for is:
\begin{multline} \label{eq:Limit_final2} 
\lim_{x\to 0}\Phi\left(e^{2\pi\ii\,x},-2k+1,n+1\right)+\frac{\pi\ii\,x}{k}\Phi\left(e^{2\pi\ii\,x},-2k,n+1\right)=\frac{(n+1)^{2k-1}}{2}-\frac{(n+1)^{2k}}{2k}\\-\frac{\ii}{2}\int_{0}^{\infty}(1-\coth{\pi x})\left((n+1+\ii\,x)^{2k-1}-(n+1-\ii\,x)^{2k-1}\right)\,dx
\end{multline}
\indent To summarize the results, the key summation has an integral representation given by:
\begin{multline} \label{eq:soma_final1} 
\sum_{j=0}^{k}\frac{(2\pi\ii\,n)^{-2j}\,\zeta(2j)}{(2k-2j)!}=\frac{1}{4n(2k-1)!}+\frac{1}{4(2k-1)!}\left(\frac{1}{n+1}-\frac{1}{k}\right)\left(1+\frac{1}{n}\right)^{2k}\\-\frac{\ii\,n^{-2k}}{4(2k-1)!}\int_{0}^{\infty}(1-\coth{\pi x})\left((n+1+\ii\,x)^{2k-1}-(n+1-\ii\,x)^{2k-1}\right)\,dx
\end{multline}

\subsection{Zeta relation to partial sums}
\indent Going back to the Faulhaber formula,
\begin{equation} \nonumber
H_{-2k+1}(n)=\frac{n^{2k-1}}{2}+2(2k-1)!(2\pi\ii)^{-2k}\zeta(2k)-2(2k-1)!\,n^{2k}\sum_{j=0}^{k}\frac{(2\pi\ii\,n)^{-2j}\zeta(2j)}{(2k-2j)!}
\end{equation}
\indent Replacing \eqrefe{soma_final1} into the above and simplifying, one has,
\begin{multline} \nonumber
H_{-2k+1}(n)=-\frac{(n+1)^{2k-1}}{2}+\frac{(n+1)^{2k}}{2k}+2(2k-1)!(2\pi\ii)^{-2k}\zeta(2k)+\\+\frac{\ii}{2}\int_{0}^{\infty}(1-\coth{\pi x})\left((n+1+\ii\,x)^{2k-1}-(n+1-\ii\,x)^{2k-1}\right)\,dx
\end{multline}
\indent As mentioned in \citena{Faulhaber}, for this formula to hold for every $k$, $\ii^{-2k}$ must be replaced with $\cos{k\pi}$ (note it is the real part of $(-1)^{-k}$). With that and also a transformation of the variable, $2k-1 \mapsto k$, one obtains:
\begin{multline} \nonumber
H_{-k}(n)=\frac{(n+1)^{k+1}}{k+1}-\frac{(n+1)^{k}}{2}+2\,k!(2\pi)^{-k-1}\cos{\frac{(k+1)\pi}{2}}\zeta(k+1)\\+\frac{\ii}{2}\int_{0}^{\infty}(1-\coth{\pi x})\left((n+1+\ii\,x)^{k}-(n+1-\ii\,x)^{k}\right)\,dx
\end{multline}
\indent Finally, since the equation below,
\begin{equation} \nonumber
2\,\Gamma(k+1)(2\pi)^{-k-1}\cos{\frac{(k+1)\pi}{2}}\zeta(k+1)=\zeta(-k) \text{,}
\end{equation}
\noindent is the Riemann functional equation, the final formula is:
\begin{multline} \label{eq:Zeta_HN_relation}
H_{-k}(n)=\frac{(n+1)^{k+1}}{k+1}-\frac{(n+1)^{k}}{2}+\zeta(-k)\\+\frac{\ii}{2}\int_{0}^{\infty}(1-\coth{\pi x})\left((n+1+\ii\,x)^{k}-(n+1-\ii\,x)^{k}\right)\,dx
\end{multline}
\indent One of the advantages of this formula, over \eqrefe{HN_final_v2} and \eqrefe{HN_final_v1}, is that it allows the partial Taylor series expansion of $H_{-k}(n)$ about $n=0$ to be obtained, even when $\Re(k)<0$. When $k$ is a positive integer, the series expansion about $n=0$ gives the Faulhaber formula, like \eqrefe{HN_final_v2} and \eqrefe{HN_final_v1}.\\

Relation \eqrefe{Zeta_HN_relation} provides the analytic continuation of $H_{-k}(n)$ to the complex plane on both parameters, $k$ and $n$:
\begin{equation} \label{eq:partial_Zeta_sums}
\sum_{q=1}^{n}q^k=\zeta(-k)-\zeta(-k,n+1)
\end{equation}
\indent One exception is $n=-1$ if $\Re(k)<0$, but there may be other possible singularities.

\subsection{$\zeta(k)$ formula}
Using equation \eqrefe{Zeta_HN_relation} with $n=0$, one obtains an integral representation for the Riemann zeta function valid in the complex plane, except its pole:
\begin{equation} \label{eq:Zeta}
\zeta(k)=\frac{1}{k-1}+\frac{1}{2}-\frac{\ii}{2}\int_{0}^{\infty}(1-\coth{\pi x})\left((1+\ii\,x)^{-k}-(1-\ii\,x)^{-k}\right)\,dx
\end{equation}
\indent The above formula can be used to easily find the derivatives of the zeta function at 0. If $q$ is a positive integer,
\begin{equation} \nonumber
\frac{\text{d}^q\,\zeta(x)}{\text{d}\,x^q}\bigg|_{x=0}=-q!-\frac{\ii(-1)^q}{2}\int_{0}^{\infty}(1-\coth{\pi x})\left(\log^q(1+\ii\,x)-\log^q(1-\ii\,x)\right)\,dx
\end{equation}

\section{Hurwitz zeta function}
Let us use \eqrefe{Zeta_HN_relation} to demonstrate how to obtain a similar  relation for the Hurwitz zeta and the generalized harmonic progressions.\\

When $v=1$, one has a special case of the identity \eqrefe{Soma_Hurwitz}, which is:
\begin{equation} \label{eq:Soma_Hurwitz_spec}
\zeta(-k,b+1)=-\frac{b^{k+1}}{k+1}+k!\,\sum_{j=0}^{k}\frac{\zeta(-j)\,b^{k-j}}{j!(k-j)!} \text{,} 
\end{equation}

\subsection{Hurwitz relation to partial sums}
First, recalling the expression in \eqrefe{Zeta_HN_relation}, and swapping $k$ for $j$:
\begin{multline} \label{eq:starting}
\sum_{q=1}^{n}q^{j}=\frac{(n+1)^{j+1}}{j+1}-\frac{(n+1)^{j}}{2}+\zeta(-j)\\+\frac{\ii}{2}\int_{0}^{\infty}(1-\coth{\pi x})\left((n+1+\ii\,x)^{j}-(n+1-\ii\,x)^{j}\right)\,dx \text{}
\end{multline}
\indent Next, the above is summed over $j$ in accordance with identity \eqrefe{Soma_Hurwitz_spec}:
\begin{multline} \nonumber
\sum_{q=1}^{n}\sum_{j=0}^{k}\frac{k!\,q^{j}\,b^{k-j}}{j!(k-j)!}=\sum_{j=0}^{k}\frac{k!\,b^{k-j}}{j!(k-j)!}\frac{(n+1)^{j+1}}{j+1}-\frac{1}{2}\sum_{j=0}^{k}\frac{k!\,(n+1)^{j}\,b^{k-j}}{j!(k-j)!}+k!\,\sum_{j=0}^{k}\frac{\zeta(-j)\,b^{k-j}}{j!(k-j)!}\\
+\frac{\ii}{2}\int_{0}^{\infty}(1-\coth{\pi x})\sum_{j=0}^{k}\frac{k!\,b^{k-j}}{j!(k-j)!}\left((n+1+\ii\,x)^{j}-(n+1-\ii\,x)^{j}\right)\,dx \text{}
\end{multline}\\
\indent Now, simplifying with the Newton's binomial,
\begin{multline} \nonumber
\sum_{q=1}^{n}(q+b)^{k}=\frac{-b^{k+1}+(n+1+b)^{k+1}}{k+1}-\frac{(n+1+b)^{k}}{2}+\frac{b^{k+1}}{k+1}+\zeta(-k,b+1)\\
+\frac{\ii}{2}\int_{0}^{\infty}(1-\coth{\pi x})\left((n+1+b+\ii\,x)^{k}-(n+1+b-\ii\,x)^{k}\right)\,dx \text{,}
\end{multline}
\noindent one finally finds that for all complex $k \neq -1$,
\begin{multline} \label{eq:Hurwitz_zeta_rel1}
\sum_{q=0}^{n}(q+b)^{k}=\frac{(n+1+b)^{k+1}}{k+1}-\frac{(n+1+b)^{k}}{2}+\zeta(-k,b)\\
+\frac{\ii}{2}\int_{0}^{\infty}(1-\coth{\pi x})\left((n+1+b+\ii\,x)^{k}-(n+1+b-\ii\,x)^{k}\right)\,dx
\end{multline}
\indent Aside from any singularities, relation \eqrefe{Hurwitz_zeta_rel1} provides the analytic continuation of the summation on the left-hand side to the complex plane on all parameters, $k$, $b$ and $n$:
\begin{equation} \label{eq:partial_Hurwitz_sums}
\sum_{q=0}^{n}(q+b)^k=\zeta(-k,b)-\zeta(-k,n+1+b)
\end{equation}

\subsection{$\zeta(-k,b)$ formula}
The easiest way to derive a formula for $\zeta(-k,b)$ is simply to set $n=-1$ in relation \eqrefe{Hurwitz_zeta_rel1}:
\begin{equation} \label{eq:Hurwitz}
\zeta(-k,b)=-\frac{b^{k+1}}{k+1}+\frac{b^{k}}{2}-\frac{\ii}{2}\int_{0}^{\infty}(1-\coth{\pi x})\left((b+\ii\,x)^{k}-(b-\ii\,x)^{k}\right)\,dx
\end{equation}
\indent This formula is valid in the complex plane (except $k=-1$, or $b=0$ if $\Re(k)<0$). It could have also been derived using equations \eqrefe{Zeta_HN_relation} and \eqrefe{partial_Zeta_sums}. 

\section{The polylogarithm}
The insight on how to derive this next formula comes from noticing the patterns in the formulae so far. When a transformation was applied to the partial sums of the zeta function, like,
\begin{equation} \nonumber
\sum_{j=0}^{k}\frac{k!\,b^{k-j}}{j!\,(k-j)!}\left(\sum_{q=1}^{n}q^{j}\right)  \text{, one obtained } \sum_{q=1}^{n}(q+b)^{k} \text{,}
\end{equation}
\noindent which are the partial sums of the Hurwitz zeta.\\

So, the next transformation needed in order to obtain the partial sums of the polylogarithm function\citesup{Abramowitz} is:
\begin{equation} \nonumber
\sum_{q=1}^{n}q^{k}e^{z\,q}=\sum_{q=1}^{n}q^{k}\sum_{j=0}^{\infty}\frac{(z\,q)^j}{j!}=\sum_{j=0}^{\infty}\frac{z^j}{j!}\sum_{q=1}^{n}q^{j+k}=\sum_{j=k}^{\infty}\frac{z^{j-k}}{(j-k)!}\sum_{q=1}^{n}q^{j} \text{}
\end{equation}
\indent This transformation must be applied to each term of \eqrefe{starting}. Let us do it by parts, but first let us get acquainted with the function and its integral. If $k$ is a non-negative integer:
\begin{equation} \nonumber
\sum_{j=k}^{\infty}\frac{x^{j}}{(j-k)!}=x^k\,e^x \text{, and } \sum_{j=k}^{\infty}\frac{x^{j+1}}{(j+1)(j-k)!}=\int_{0}^{x}v^k\,e^v\,dv=(-1)^k\left(\Gamma(k+1,-x)-k!\right)
\end{equation}
\indent Therefore, the first term is:
\begin{equation} \nonumber
\sum_{j=k}^{\infty}\frac{z^{j-k}}{(j-k)!}\frac{(n+1)^{j+1}}{j+1}=-(-z)^{-k-1}\left(\Gamma(k+1,-z(n+1))-k!\right)
\end{equation}
\indent And the second term is:
\begin{equation} \nonumber
-\frac{1}{2}\sum_{j=k}^{\infty}\frac{z^{j-k}(n+1)^j}{(j-k)!}=-\frac{(n+1)^k e^{z(n+1)}}{2}
\end{equation}
\indent The next term warrants an entire section.

\subsection{Zeta at the negative integers}
For $k$ a non-negative integer, a closed form for the following power series is needed:
\begin{equation} \nonumber
\sum_{j=k}^{\infty}\frac{z^{j-k}}{(j-k)!}\zeta(-j)
\end{equation}
\indent One can, therefore, start from the generating function of the zeta function at the negative integers, whose $k$-th derivative gives the above:
\begin{equation} \nonumber
\sum_{j=0}^{\infty}\frac{x^{j}}{j!}\zeta(-j)=-\frac{1}{2}+\sum_{j=1}^{\infty}\frac{x^{2j-1}}{(2j-1)!}\zeta(-2j+1)=-\frac{1}{2}-\sum_{j=1}^{\infty}\frac{B_{2j}\,x^{2j-1}}{(2j)!}=-\frac{1}{2}+\frac{1}{x}-\frac{1}{2}\coth{\frac{x}{2}}
\end{equation}\\
\indent The derivatives of the hyperbolic cotangent can be calculated with the identities from section \secrefe{Trigon}:
\begin{equation} \nonumber
\frac{\text{d}^k}{\text{d}\,x^k}\left(-\frac{1}{2}+\frac{1}{x}-\frac{1}{2}\coth{\frac{x}{2}}\right)\bigg|_{x=z}=(-1)^{k}k!\,z^{-k-1}-\delta_{0\,k}-(-1)^k\,\mathrm{Li}_{-k}\left(e^{-z}\right)
\end{equation}\\
\indent The Kronecker delta is a problem, but fortunately it is resolved by the following equivalence:
\begin{equation} \nonumber
-\delta_{0\,k}-(-1)^k\,\mathrm{Li}_{-k}\left(e^{-z}\right)=\mathrm{Li}_{-k}\left(e^{z}\right)
\end{equation}

\subsection{Polylogarithm relation to partial sums}
The last term is simple:
\begin{multline} \nonumber
\frac{\ii}{2}\int_{0}^{\infty}(1-\coth{\pi x})\sum_{j=k}^{\infty}\frac{z^{j-k}}{(j-k)!}\left((n+1+\ii\,x)^{j}-(n+1-\ii\,x)^{j}\right)\,dx\\ 
=\frac{\ii\,e^{z(n+1)}}{2}\int_{0}^{\infty}(1-\coth{\pi x})\left(e^{z\,\ii\,x}(n+1+\ii\,x)^{k}-e^{-z\,\ii\,x}(n+1-\ii\,x)^{k}\right)\,dx
\text{}
\end{multline}
\indent When everything is put together, one obtains:
\begin{multline} \label{eq:Polylog_rel_final}
\sum_{q=1}^{n}q^{k}e^{z\,q}=-\frac{(n+1)^k e^{z(n+1)}}{2}-(-z)^{-k-1}\Gamma(k+1,-z(n+1))+\mathrm{Li}_{-k}\left(e^{z}\right)\\ 
+\frac{\ii\,e^{z(n+1)}}{2}\int_{0}^{\infty}(1-\coth{\pi x})\left(e^{z\,\ii\,x}(n+1+\ii\,x)^{k}-e^{-z\,\ii\,x}(n+1-\ii\,x)^{k}\right)\,dx
\text{}
\end{multline}
\indent This relation provides the analytic continuation of the summation on the left-hand side to the complex plane on all parameters, $k$, $z$ and $n$:
\begin{equation} \label{eq:partial_polylog_sum}
\sum_{q=1}^{n}q^{k}e^{z\,q}=\mathrm{Li}_{-k}\left(e^{z}\right)-e^{z(n+1)}\Phi\left(e^{z},-k,n+1\right)
\end{equation}

\subsection{$\mathrm{Li}_{-k}\left(e^{z}\right)$ formula}
The simplest formula for $\mathrm{Li}_{-k}\left(e^{z}\right)$ is obtained by setting $n=0$ in relation \eqrefe{Polylog_rel_final},
\begin{multline} \label{eq:Polylog_final}
\mathrm{Li}_{-k}\left(e^{z}\right)=\frac{e^{z}}{2}+(-z)^{-k-1}\Gamma(k+1,-z)\\ 
-\frac{\ii\,e^{z}}{2}\int_{0}^{\infty}(1-\coth{\pi x})\left(e^{z\,\ii\,x}(1+\ii\,x)^{k}-e^{-z\,\ii\,x}(1-\ii\,x)^{k}\right)\,dx
\text{,}
\end{multline}
\noindent which should be valid in the complex plane, except when $e^{z}=1$.


\section{Lerch $\Phi$ function}
For the Lerch $\Phi$ function\citesup{Lerch1887} the process is the very same, but this time the identity \eqrefe{Soma_Polylog} is used, instead of \eqrefe{Soma_Hurwitz}.

\subsection{Lerch $\Phi$ relation to partial sums}
Summing both sides of equation \eqrefe{Polylog_rel_final} (with $k$ replaced by $j$), over $j$ per the transformation \eqrefe{Soma_Polylog},
\begin{multline} \nonumber
\sum_{q=1}^{n}e^{z\,q}\sum_{j=0}^{k}\frac{k!\,q^{j}\,b^{k-j}}{j!(k-j)!}=-\frac{e^{z(n+1)}}{2}\sum_{j=0}^{k}\frac{k!\,(n+1)^j\,b^{k-j}}{j!(k-j)!}
\\-\sum_{j=0}^{k}\frac{k!\,b^{k-j}}{j!(k-j)!}(-z)^{-j-1}\Gamma(j+1,-z(n+1))+\sum_{j=0}^{k}\frac{k!\,\mathrm{Li}_{-j}\left(e^{z}\right)\,b^{k-j}}{j!(k-j)!}\\ 
+\frac{\ii\,e^{z(n+1)}}{2}\int_{0}^{\infty}(1-\coth{\pi x})\sum_{j=0}^{k}\frac{k!\,b^{k-j}}{j!(k-j)!}\left(e^{z\,\ii\,x}(n+1+\ii\,x)^{j}-e^{-z\,\ii\,x}(n+1-\ii\,x)^{j}\right)\,dx
\text{,}
\end{multline}
\noindent one concludes that:
\begin{multline} \label{eq:Lerch_rel_final}
\sum_{q=0}^{n}(q+b)^{k}\,e^{z\,q}=-\frac{(n+1+b)^k e^{z(n+1)}}{2}-(-z)^{-k-1}e^{-z\,b}\,\Gamma(k+1,-z(n+1+b))+\Phi(e^{z},-k,b)\\ 
+\frac{\ii\,e^{z(n+1)}}{2}\int_{0}^{\infty}(1-\coth{\pi x})\left(e^{z\,\ii\,x}(n+1+b+\ii\,x)^{k}-e^{-z\,\ii\,x}(n+1+b-\ii\,x)^{k}\right)\,dx
\text{}
\end{multline}
\indent Aside from singularities, relation \eqrefe{Lerch_rel_final} provides the analytic continuation of the summation on the left-hand side to the complex plane on all parameters, $k$, $z$, $b$ and $n$:
\begin{equation} \label{eq:partial_Lerch_sums}
\sum_{q=0}^{n}(q+b)^k\,e^{z\,q}=\Phi(e^{z},-k,b)-e^{z(n+1)}\,\Phi(e^{z},-k,n+1+b)
\end{equation}

\subsection{$\Phi(e^{z},-k,b)$ formula}
Setting $n=-1$ in relation \eqrefe{Lerch_rel_final} is the easiest way to derive a Lerch $\Phi$ formula valid in the complex plane (except for occasional singularities):
\begin{multline} \label{eq:Lerch_final}
\Phi(e^{z},-k,b)=\frac{b^k}{2}+(-z)^{-k-1}e^{-z\,b}\,\Gamma(k+1,-z\,b)\\ 
-\frac{\ii}{2}\int_{0}^{\infty}(1-\coth{\pi x})\left(e^{z\,\ii\,x}(b+\ii\,x)^{k}-e^{-z\,\ii\,x}(b-\ii\,x)^{k}\right)\,dx
\text{}
\end{multline}

\section{Different formulae}
These formulae must not hold everywhere (except for the parameter $k$). The alternative Lerch $\Phi$ formula used in section \secrefe{Lerch_diff_lim_sec} came from this method. If the steps outlined previously are repeated using the expression that was created in \eqrefe{HN_final_v2}, the following formulae are obtained.

\subsection{Hurwitz zeta function}
After all is put together, one finds that for all complex $k \neq -1$:
\begin{multline} \label{eq:Hurwitz_zeta_rel2} \nonumber
\sum_{q=0}^{n}(q+b)^{k}=\frac{(n+b)^{k+1}}{k+1}+\frac{(n+b)^{k}}{2}+\zeta(-k,b)\\
-n\int_{0}^{\pi/2}\frac{1-\coth{\left(\pi\,n\tan{v}\right)}}{(\cos{v})^{2}}\left((n+b)^2+(n\tan{v})^2\right)^{k/2}\sin{\left(k\arctan{\frac{n\tan{v}}{n+b}}\right)}\,dv \text{,}
\end{multline}
\indent It is not possible to turn this integral into the one from equation \eqrefe{Hurwitz_zeta_rel1}, and the reason is because two of the terms outside of the integral differ between the two formulae.

\subsection{The Polylogarithm}
Using equation \eqrefe{HN_final_v2}, this is the polylogarithm formula that one ends up with:
\begin{multline} \nonumber
\sum_{q=1}^{n}q^{k}e^{z\,q}=\frac{n^k\,e^{z\,n}}{2}-(-z)^{-k-1}\Gamma(k+1,-z\,n)+\mathrm{Li}_{-k}\left(e^{z}\right)\\ 
-n^{k+1}e^{z\,n}\int_{0}^{\pi/2}(1-\coth{(\pi\,n\tan{v})})\frac{\sin{\left(k\,v+z\,n\tan{v}\right)}}{(\cos{v})^{k+2}}\,dv
\text{}
\end{multline}

\subsection{Lerch $\Phi$ function} \label{Lerch_used}
Using the polylogarithm relation from the previous subsection one obtains:
\begin{multline} \nonumber
\sum_{q=0}^{n}(q+b)^{k}e^{z\,q}=\frac{(n+b)^k\,e^{z\,n}}{2}-(-z)^{-k-1}e^{-z\,b}\,\Gamma(k+1,-z(n+b))+\Phi\left(e^{z},-k,b\right)\\ 
-n\,e^{z\,n}\int_{0}^{\pi/2}\frac{1-\coth{(\pi\,n\tan{v})}}{(\cos{v})^{2}}\left((n+b)^2+(n\,\tan{v})^2\right)^{k/2}\sin{\left(k\arctan{\frac{n\tan{v}}{n+b}}+z\,n\tan{v}\right)}\,dv
\text{}
\end{multline}

\end{document}